\documentclass[12pt,reqno]{amsart}

\usepackage{enumitem}   

\usepackage[colorlinks = true,
            linkcolor = blue,
            urlcolor  = blue,
            citecolor = blue,
            anchorcolor = blue]{hyperref}
						
\newtheorem{theorem}{Theorem}

{Proposition}
\newtheorem{corollary}
{Corollary}

\newtheorem{lemma}{Lemma}

\newfont{\bb}{msbm10 at 12pt}

\def\<{\langle}     
\def\>{\rangle}

\newcommand{\bal}{\begin{align}}      \newcommand{\eal}{\end{align}}
\newcommand{\ba}{\begin{array}}      \newcommand{\ea}{\end{array}}
\newcommand{\bc}{\begin{center}}     \newcommand{\ec}{\end{center}}
\newcommand{\be}{\begin{enumerate}}  \newcommand{\ee}{\end{enumerate}}
\newcommand{\beq}{\begin{eqnarray}}  \newcommand{\eeq}{\end{eqnarray}}
\newcommand{\beQ}{\begin{eqnarray*}} \newcommand{\eeQ}{\end{eqnarray*}}
\newcommand{\bi}{\begin{itemize}}    \newcommand{\ei}{\end{itemize}}
\newcommand{\bt}{\begin{tabular}}    \newcommand{\et}{\end{tabular}}
\newcommand{\bdm}{\begin{displaymath}} \newcommand{\edm}{\end{displaymath}}

\newcommand{\Ss}{R\!\!\!\!/\,}

\begin{document}

\title[]{On the ADM mass of critical area-normalized capacitors}       

\author{Simon Raulot}
\address[Simon Raulot]{Univ Rouen Normandie, CNRS, Normandie Univ, LMRS UMR 6085, F-76000 Rouen, France}
\email{simon.raulot@univ-rouen.fr}

\begin{abstract}
In this note, we prove mass-capacity inequalities for asymptotically flat manifolds whose boundary capacity potential satisfies an overdetermined problem, referred to as critical area-normalized capacitors. As a consequence, we obtain uniqueness results for the Schwarzschild metric,  from which improvements in the uniqueness theorems for spin asymptotically flat spacetimes containing a connected photon surface, as well as for spin asymptotically flat static manifolds with boundary are obtained. 
\end{abstract}

\keywords{}

\subjclass{53C20, 52C24, 53C27, 83C57}

\thanks{}

\date{\today}   

\maketitle 
\pagenumbering{arabic}


\section{Statements of the main results}


A smooth, connected and oriented $n$-dimensional Riemannian manifold $(M^n,g)$ is said to be asymptotically flat (of order $\tau$) if there exists a compact subset $K\subset M$ such that $M\setminus K$ is a finite disjoint union of ends $M_k$, each of them being diffeomorphic to $\mathbb{R}^n$ minus a closed ball $\overline{B}$ by a coordinate chart in which the components of the metric satisfy 
\begin{eqnarray*}
g_{ij}=\delta_{ij}+O_2\big(r^{-\tau}\big)
\end{eqnarray*}
for $i,j=1,...,n$, $\tau>(n-2)/2$ and the scalar curvature $R\in L^1(M)$. Here $O_2(r^{-\tau})$ refers to a real-valued function $f$ such that
\begin{eqnarray*}
|f(x)|+r|\partial f(x)|+r^2|\partial^2f(x)|\leq Cr^{-\tau}
\end{eqnarray*}
as $r$ goes to infinity, for some constant $C>0$ and where $\partial$ is the standard derivative in the Euclidean space. Such a coordinate chart is often referred to as a chart at infinity. In the following, we assume that $k=1$ and the general case can be treated in a similar way.  

On an asymptotically flat manifold $(M^n,g)$, the {\it ADM mass} is defined by 
\begin{eqnarray*}
m:=\frac{1}{2(n-1) \omega_{n-1}}\lim_{r\rightarrow+\infty}\sum_{i,j=1}^n\int_{\mathcal{S}_r}(g_{ij,i}-g_{ii,j}) \frac{x^j}{r}d\overline{\sigma}_{r}
\end{eqnarray*}
where $\mathcal{S}_r$ stands for a coordinate sphere of radius $r>0$, $d\overline{\sigma}_{r}$ its Euclidean Riemannian volume element, $g_{ij,s}$ the derivative of the metric components in the coordinate chart and $\omega_{n-1}$ is the volume of the $(n-1)$-dimensional sphere. Although this definition seems to depend on a specific choice of a coordinate chart, it is, as independently proven by Bartnik \cite{Bartnik1} and Chruściel \cite{Chrusciel1}, a well-defined geometric invariant. The Riemannian positive mass theorem by Schoen and Yau \cite{SchoenYau1,SchoenYau2} for $3\leq n\leq 7$ (but see \cite{Lohkamp2,SchoenYau3} for the general case) and by Witten \cite{Witten1} for spin manifolds asserts that if the scalar curvature of $(M^n,g)$ is nonnegative, then $m\geq 0$ and it is zero if and only if $(M^n,g)$ is isometric to the Euclidean space. 

If now we assume in addition that $M$ has a compact inner boundary $\Sigma$ with nonpositive mean curvature, one can show that the ADM mass has to be positive. In this situation, Bray \cite{Bray} (for $n=3$) and Hirsch-Miao \cite{HirschMiao1} (for manifolds where the positive mass theorem holds) proved that $m\geq \mathcal{C}(\Sigma,M)$ with equality if, and only if, $(M^n,g)$ is isometric to the exterior of the rotationally symmetric minimal sphere in the Riemannian Schwarzschild manifold of corresponding mass. Here $\mathcal{C}(\Sigma,M)$ is the {\it boundary capacity} of $\Sigma$ in $(M^n,g)$ defined as
\begin{eqnarray*}\label{DefCap}
\mathcal{C}(\Sigma,M):=\inf_{f\in\mathcal{M}_{1,0}}\left\{\frac{1}{(n-2)\omega_{n-1}}\int_M|\nabla f|^2d\mu\right\}
\end{eqnarray*}
where $\mathcal{M}_{1,0}$ denotes the set of all smooth functions on $M$ which are exactly $1$ on $\Sigma$ and approach $0$ towards infinity in the end and $\nabla$ is the gradient of $(M^n,g)$. This minimum is achieved by a smooth harmonic function $\Phi\in\mathcal{M}_{1,0}$ usually called the {\it boundary capacity potential} which therefore satisfies the following boundary value problem  
$$
\left\lbrace
\begin{array}{ll}
\Delta\Phi = 0 & \text{ on } M\\
\Phi=1 & \text{ at } \Sigma\\
\Phi(x)\rightarrow 0 & \text{ as }|x|\rightarrow\infty.
\end{array}
\right.
$$
Here $\Delta$ is the Laplace operator on $(M^n,g)$. Note that since $\Phi$ is harmonic, it holds that 
\begin{eqnarray}\label{ExpansionBCP}
\Phi(x)=\frac{\mathcal{C}(\Sigma,M)}{r^{n-2}}+O_2(r^{-(n-2+\varepsilon)})
\end{eqnarray}
as $r\rightarrow\infty$ and where $\varepsilon>0$ (see \cite[Theorem 2.2]{AgostinianiMazzieriOronzio} for example) from which we deduce that 
\begin{eqnarray}\label{DefCap}
\mathcal{C}(\Sigma,M)=-\frac{1}{(n-2)\omega_{n-1}}\int_\Sigma\frac{\partial\Phi}{\partial\nu}d\sigma
\end{eqnarray}
where $\nu$ is the unit normal pointing inside $M$. 

In this note, we will study the mass of asymptotically flat manifolds for which the boundary capacity potential satisfies an additional boundary condition. More precisely, we will assume that there exists a constant $\Lambda\in\mathbb{R}$ such that
\begin{equation}\label{AdditionalBC}
{\frac{\partial \Phi}{\partial\nu}}\Big|_{\Sigma}=-\frac{1}{2}\frac{n-2}{n-1}\Lambda H
\end{equation}
where $H$ denotes the mean curvature of $\Sigma$ in $(M^n,g)$. Here the mean curvature is the component of the mean curvature vector field with respect to $-\nu$. Such a manifold will be referred to as a {\it critical area-normalized capacitor} since as noticed in \cite{FallMinlendRatzkin}, they arise as critical points of a certain area-normalized domain functional. It is immediate from (\ref{DefCap}) that the capacity of such critical points is related to the total mean curvature of $\Sigma$ in $M$ by the formula
\begin{eqnarray*}
\mathcal{C}(\Sigma,M)=\frac{\Lambda}{2(n-1)\omega_{n-1}}\int_\Sigma H\,d\sigma.
\end{eqnarray*}

The prototypical examples of critical area-normalized capacitors are given by the Riemannian Schwarzschild manifolds. These are the one-parameter family of manifolds defined for $m\in\mathbb{R}$ by 
$$
\mathbb{M}^n_m:=
\left\lbrace
\begin{array}{ll}
\mathbb{R}^n\setminus\{0\} & \text{ if } m>0\\
\mathbb{R}^n & \text{ if } m=0\\
\mathbb{R}^{n}\setminus\big\{r\leq\big(|m|/2\big)^{1/(n-2)}\big\} & \text{ if } m<0
\end{array}
\right.
$$
and equipped with the Riemannian metric 
\begin{eqnarray*}\label{SchwarzschildMetric}
g_m=\Big(1+\frac{m}{2r^{n-2}}\Big)^{\frac{4}{n-2}}\delta
\end{eqnarray*}
where $r:=|x|$ is the Euclidean radius for $x\in\mathbb{M}^n_m$ and $\delta$ is the Euclidean metric. For $r_0\in(r_\ast,\infty)$ with $r_\ast=0$ if $m\geq 0$ and $r_\ast=(|m|/2)^{1/(n-2)}$ if $m<0$, we consider the exterior of the region outside a rotationally symmetric sphere defined by
\begin{eqnarray*}
\mathbb{M}_m^n(r_0):=\Big\{x\in\mathbb{M}^n_m\,/\,r\geq r_0\Big\}.
\end{eqnarray*}
This is an $n$-dimensional, spin, complete and asymptotically flat manifold with zero scalar curvature. Moreover, its inner boundary $\Sigma_{r_0}$ with induced metric $\gamma_{r_0}$ is isometric to a round sphere with radius
\begin{eqnarray}\label{RadiusSS}
r_{m,r_0}:=r_0\big(1+\frac{m}{2r_0^{n-2}}\big)^{\frac{2}{n-2}}
\end{eqnarray}
and constant mean curvature 
\begin{eqnarray}\label{MeanCurvatureSS}
H_{m,r_0}=\frac{n-1}{r_0}\big(1-\frac{m}{2r_0^{n-2}}\big)\big(1+\frac{m}{2r_0^{n-2}}\big)^{-\frac{n}{n-2}}.
\end{eqnarray} 
On the other hand, it is not difficult to check that the boundary capacity potential of $\Sigma_{r_0}$ in $\big(\mathbb{M}^n_m,g_m\big)$ is given by 
\begin{eqnarray}\label{CapacityPotentialSS}
\Phi_{r_0}(x)=\big(1+\frac{m}{2r_0^{n-2}}\big)\big(1+\frac{m}{2r^{n-2}}\big)^{-1}\Big(\frac{r_0}{r}\Big)^{n-2}
\end{eqnarray}
and so it holds on $\Sigma_{r_0}$ that 
\begin{eqnarray*}\label{NormalDerivativeBCP}
\frac{\partial\Phi_{r_0}}{\partial\nu}=-\frac{n-2}{r_0}\big(1+\frac{m}{2r_0^{n-2}}\big)^{-\frac{n}{n-2}}.
\end{eqnarray*}
In particular we have
\begin{eqnarray}\label{SchwarzschildCapacitors}
\frac{\partial\Phi_{r_0}}{\partial\nu}=-\frac{1}{2}\frac{n-2}{n-1}\Lambda_{m,r_0}H_{m,r_0}\quad\text{with}\quad\Lambda_{m,r_0}=\frac{2}{1-\frac{m}{2r^{n-2}_0}}
\end{eqnarray}
so that the Riemannian Schwarzschild manifold is actually a critical area-normalized capacitor for all $r_0>r_\ast$ if $m<0$ and for $r_0\neq (m/2)^{1/(n-2)}$ when $m\geq 0$. 

It appears that, under some natural geometrical assumptions, the mass of a critical area-normalized capacitor satisfies an explicit lower bound. In the following, $\Ss$ denotes the scalar curvature of $(\Sigma^{n-1},\gamma)$ and $\gamma=g_{|\Sigma}$. 
\begin{theorem}\label{SpinCase}
Let $(M^n,g)$ be a spin critical area-normalized capacitor with a connected boundary $\Sigma^{n-1}:=\partial M^n$. Assume that $R\geq 0$ as well as $H>0$ and that there exists $c>1$ such that 
\begin{equation}\label{ScalarMeanHypothese}
\inf_\Sigma(\Ss)\geq\frac{n-2}{n-1}c\max_{\Sigma}(H^2).
\end{equation}
Then 
\begin{equation}\label{MassCapIneq}
m\geq\frac{1}{1+\frac{\Lambda}{c-1}}\,\mathcal{C}(\Sigma,M)
\end{equation}
and equality occurs if and only if $(M^n,g)$ is isometric to the region outside a rotationally symmetric sphere in an $n$-dimensional Schwarzschild manifold. 
\end{theorem}

Remark that the assumption (\ref{ScalarMeanHypothese}) appears in \cite{McCormickMiao} where a Penrose-like inequality is proved. When $n=3$, we can relax the assumption (\ref{ScalarMeanHypothese}). 
\begin{theorem}\label{ThreeDimensionalCase}
Let $(M^n,g)$ be a critical area-normalized capacitor with a connected boundary $\Sigma^{n-1}:=\partial M^n$. Assume that $R\geq 0$ as well as $H>0$ and that there exists $c>1$ such that 
\begin{equation}\label{ScalarMeanHypothese1}
\Ss\geq\frac{n-2}{n-1}cH^2.
\end{equation}
If one of the following conditions holds:
\begin{itemize}
\item $n=3$;
\item $(\Sigma^{n-1},\gamma)$ admits an isometric embedding in the flat $n$-dimensional Euclidean space;
\end{itemize} 
then the inequality (\ref{MassCapIneq}) holds and equality occurs if and only if $(M^n,g)$ is isometric to the region outside a rotationally symmetric sphere exterior in an $n$-dimensional Schwarzschild manifold.
\end{theorem}

We refer to Section \ref{TheProofs} for the proofs. As we will see in the next section, these results apply to get rigidity results for the Schwarzschild manifold when considering static vacuums. For this, we remark that it is straightforward to check that $(M^n,g)$ is a critical area-normalized capacitor if and only if for $0<\alpha<1$, the smooth function
\begin{eqnarray*}
V:=1+(\alpha-1)\Phi:M\longrightarrow[\alpha,1)
\end{eqnarray*}
satisfies the overdetermined problem 
\begin{equation}\label{OverBP}
\Delta V=0\text{ on } M,\quad V_{|\Sigma}=\alpha,\quad V(x)\underset{|x|\rightarrow\infty}{\rightarrow}1,\quad\frac{\partial V}{\partial\nu}\Big|_{\Sigma}=
\frac{1}{2}\frac{n-2}{n-1}\Gamma H V_{|\Sigma}
\end{equation}
where $\Gamma=\Lambda(\alpha^{-1}-1)$. With this alternative characterization of critical capacitors, it is worth pointing out the following equivalent formulation of Theorem \ref{SpinCase}.  
\begin{corollary}\label{MassCapacityAlternative}
Let $(M^n,g)$ be a spin asymptotically flat manifold with a compact and connecter inner boundary. Assume that $R\geq 0$ as well as $H>0$ and that (\ref{ScalarMeanHypothese}) holds. If there exists $\alpha\in(0,1)$ such that the overdetermined problem (\ref{OverBP}) admits a non-trivial solution $V\in C^\infty(M)$, then the ADM mass $m$ of $(M^n,g)$ satisfies 
\begin{equation}\label{MassCapIneq1}
m\geq\frac{1}{1+\frac{\alpha\Gamma}{(1-\alpha)(c-1)}}\,\mathcal{C}(\Sigma,M).
\end{equation}
Moreover, equality occurs if and only if $(M^n,g)$ is isometric to the region outside a rotationally symmetric sphere in an $n$-dimensional Schwarzschild manifold.
\end{corollary}

For Theorem \ref{ThreeDimensionalCase}, we get:
\begin{corollary}\label{MassCapacityAlternative1}
Let $(M^n,g)$ be an asymptotically flat manifold with a compact and connecter inner boundary. Assume that $R\geq 0$ as well as $H>0$ and that (\ref{ScalarMeanHypothese1}) holds. If one of the following conditions is fulfilled:
\begin{itemize}
\item $n=3$;
\item $(\Sigma^{n-1},\gamma)$ admits an isometric embedding in the flat $n$-dimensional Euclidean space. 
\end{itemize} 
and if there exists $\alpha\in(0,1)$ such that the overdetermined problem (\ref{OverBP}) admits a non-trivial solution $V\in C^\infty(M)$, then the ADM mass $m$ of $(M^n,g)$ satisfies the mass-capacity inequality (\ref{MassCapIneq1}). Moreover, equality occurs if and only if $(M^n,g)$ is isometric to the region outside a rotationally symmetric sphere in an $n$-dimensional Schwarzschild manifold.
\end{corollary}


\section{Rigidity results for static asymptotically flat vacuums with equipotential boundaries}


In this section, we will obtain rigidity results for the Schwarzschild manifold in the context of asymptotically flat static vacuums. Recall that a Riemannian manifold $(M^n,g)$ is a {\it static vacuum} (with zero cosmological constant) if it admits a solution $V:M\rightarrow\mathbb{R}_+$ to the static equations
\begin{eqnarray}\label{StaticEquations}
\nabla^2 V  =  V {\rm Ric}\quad\text{and}\quad\Delta V  =  0
\end{eqnarray}
where $\nabla^2$ and ${\rm Ric}$ are respectively the Hessian and the Ricci tensor of $(M^n,g)$. The function $V$ is usually refer to as a {\it static potential}. These manifolds are of particular interest in general relativity since $(M^n,g)$ is a static vacuum with static potential $V$ if, and only if, the Lorentzian manifold
\begin{eqnarray*}
\left(\mathfrak{L}^{n+1}:=\mathbb{R}\times \left(M^n\setminus V^{-1}(\{0\})\right),\mathfrak{g}:=-V^2dt^2+g\right)
\end{eqnarray*}
is a solution of the vacuum Einstein field equations with zero cosmological constant $\mathfrak{Ric}=0$ where $\mathfrak{Ric}$ denotes the Ricci tensor of $\mathfrak{g}$. When $(M^n,g)$ has a boundary $\Sigma$, this boundary is said to be {\it equipotential} if the static potential is constant along $\Sigma$. Furthermore, if $(M^n,g)$ is asymptotically flat we will always assume that the static potential $V$ is bounded so that from \cite[Proposition B.4]{HuangMartinMiao} we can assume without loss of generality that 
\begin{eqnarray*}
\lim_{|x|\rightarrow\infty} V(x)=1
\end{eqnarray*}
within a chart at infinity. We refer to the triple $(M^n,g,V)$ as an {\it asymptotically flat static vacuum with equipotential boundary}. The exterior part $(\mathbb{M}^n_m(r_0),g_m,V_m)$ of the Schwarzschild manifold of mass $m>0$ with $r_0\geq r_m:=(m/2)^{1/(n-2)}$ and 
\begin{eqnarray*}
V_m(r)=\frac{1-\frac{m}{2r^{n-2}}}{1+\frac{m}{2r^{n-2}}}
\end{eqnarray*}
provides a $1$-parameter family of examples of asymptotically flat static vacuums with equipotential boundaries. 

As a direct application of Theorem \ref{MassCapacityAlternative}, we obtain the following rigidity result for asymptotically flat static vacuums. 
\begin{theorem}\label{RigidityDim3}
Let $(M^n,g,V)$ be an asymptotically flat static vacuum with a compact and connected equipotential boundary $\Sigma^{n-1}$ and ADM mass $m$. Assume that $\Sigma$ is totally umbilical with $H>0$ and 
\begin{eqnarray}\label{ScalarMeanCurvatureConditionRig}
\Ss=\frac{n-2}{n-1}cH^2
\end{eqnarray}
for a real $c>1$. If one of the following conditions holds:
\begin{itemize}
\item $n=3$;
\item $(\Sigma^{n-1},\gamma)$ admits an isometric embedding in the flat $n$-dimensional Euclidean space;
\end{itemize} 
then $(M^n,g)$ is isometric to the region outside a rotationally symmetric sphere in the $n$-dimensional Schwarzschild manifold of mass $m>0$ and $V$ coincides with $V_m$ (up to the isometry).
\end{theorem}
Here we do not assume that the scalar curvature or the mean curvature are constant. When $n\geq 4$, one can replace the second condition in the previous theorem as follows. 
\begin{theorem}\label{RigidityGeneral}
Let $(M^n,g,V)$ be a spin asymptotically flat static vacuum with a compact and connected equipotential boundary $\Sigma^{n-1}$ and ADM mass $m$. If $\Sigma$ is totally umbilical with constant mean curvature $H>0$ satisfying (\ref{ScalarMeanCurvatureConditionRig}) for $c>1$ then $(M^n,g)$ is isometric to the region outside a rotationally symmetric sphere in the $n$-dimensional Schwarzschild manifold of mass $m>0$ and $V$ coincides with $V_m$ (up to the isometry).
\end{theorem}
For $n=3$, this directly follows from Theorem \ref{RigidityDim3}. A simple but important ingredient in the proof of these theorems is the following classical formula regarding equipotential hypersurfaces in static vacuums. 
\begin{lemma}\label{EPS_Lemma}
The static potential $V$ on an equipotential hypersurface $\Sigma$ with $V\neq 0$ on $\Sigma$ and with $H>0$ in a static vacuum satisfies
\begin{eqnarray*}
\frac{\partial V}{\partial\nu}\Big|_{\Sigma}=\frac{1}{2}\left(\frac{n-2}{n-1}(f-1)H+H^{-1}|\mathcal{O}|^2\right)V_{|\Sigma}
\end{eqnarray*}
where $f:=\frac{n-1}{n-2}\Ss/H^2$ and $\mathcal{O}$ is the trace free part of the second fundamental form of $\Sigma$ in $M$. 
\end{lemma}
{\it Proof.} Since $V$ is harmonic and $\Sigma$ is equipotential, we compute that
\begin{eqnarray*}
H\frac{\partial V}{\partial\nu}\Big|_{\Sigma}=-\nabla^2V(\nu,\nu).
\end{eqnarray*}
The expected equality follows from the first equation in (\ref{StaticEquations}) and the Gauss equation which reads 
\begin{eqnarray*}
{\rm Ric}(\nu,\nu)=\frac{1}{2}\left(\Ss-\frac{n-2}{n-2}H^2+|\mathcal{O}|^2\right),
\end{eqnarray*}
since $g$ is scalar flat. 
\hfill$\square$

\vspace{0.2cm}

{\it Proof of Theorems \ref{RigidityDim3} and \ref{RigidityGeneral}.} First note that since $\Sigma$ is equipotential with positive mean curvature and $V$ satisfies the first equation in (\ref{StaticEquations}), the Hopf boundary point lemma ensures that $V_{|\Sigma}:=\alpha>0$ is a positive constant on $\Sigma$. Then since $\Sigma$ is totally umbilical and satisfies (\ref{ScalarMeanCurvatureConditionRig}) it follows from Lemma \ref{EPS_Lemma} that 
\begin{eqnarray}\label{NormalSP}
	\frac{\partial V}{\partial\nu}\Big|_{\Sigma}=\frac{1}{2}\frac{n-2}{n-1}(c-1)HV_{|\Sigma}.
\end{eqnarray} 
On the other hand, it is a well-known fact (see \cite[p.199-200]{Lee} for example) that under our assumptions, the static potential $V$ has the following asymptotic expansion at infinity
\begin{eqnarray}\label{AsymptExpansionSP}
V(x)=1-m r^{2-n}+O_2(r^{2-n-\varepsilon}) 
\end{eqnarray}
for some $\varepsilon>0$ and where $r=|x|$ in a chart at infinity. 
Harmonicity of $V$ and the expansion (\ref{AsymptExpansionSP}) lead to the {\it Smarr formula} 
\begin{eqnarray}\label{NormalSP1}
\int_\Sigma\frac{\partial V}{\partial\nu}d\sigma=\lim_{r\rightarrow\infty}\int_{\mathcal{S}_r}\frac{\partial V}{\partial r}d\overline{\sigma}_{r}=(n-2)\omega_{n-1}m
\end{eqnarray}
which, when combined with (\ref{NormalSP}), gives that 
\begin{eqnarray*}
m=\frac{\alpha}{2}\frac{c-1}{(n-1)\omega_{n-1}}\int_\Sigma Hd\sigma>0
\end{eqnarray*}
since $\alpha>0$, $c>1$ and $H>0$. It therefore follows from this fact, the expansion (\ref{AsymptExpansionSP}) and the maximum principle that $0<\alpha<1$. This ensures that the assumptions of Corollaries \ref{MassCapacityAlternative} and \ref{MassCapacityAlternative1} are satisfied so that they apply for $\Gamma=c-1$ leading to the mass-capacity inequality $m\geq(1-\alpha)C(\Sigma,M)$. It remains to show that this is in fact an equality. This is a direct consequence of the fact that in our situation, the boundary capacity potential is given by $\Phi=(1-V)/(1-\alpha)$ so that combining formula (\ref{DefCap}) with (\ref{NormalSP}) first and then with (\ref{NormalSP1}) allows to conclude. The rigidity part of the aforementioned corollaries finish the proof.    
\hfill$\square$

\vspace{0.2cm}

To conclude this section, let us mention two direct applications of the above theorems.

\subsection{An uniqueness theorem for asymptotically flat static vacuum spacetimes with a connected photon surface as inner boundary} A direct consequence of Theorem \ref{RigidityGeneral} is the following version of the uniqueness theorem for asymptotically flat static vacuum spacetime.
\begin{theorem}\label{PhotonSurfacesUniqueness}
Let $(\mathfrak{L}^{n+1},\mathfrak{g})$ be a spin asymptotically flat static vacuum spacetime with ADM mass $m$. Assume that $(\mathfrak{L}^{n+1},\mathfrak{g})$ is geodesically complete up to its inner boundary $\partial\mathfrak{L}$ which is assumed to be a connected photon surface $P^n:=\partial\mathfrak{L}$. Assume, in addition, that $P^n$ is equipotential, outward directed and has compact time slices. Then $(\mathfrak{L}^{n+1},\mathfrak{g})$ is isometric to a suitable piece of the Schwarzschild spacetime of positive mass $m$. 
\end{theorem}
We refer to \cite{CederbaumGalloway2} for the precise definitions of the various objects involved in this statement. The proof relies on the fact that 
such a spacetime naturally leads to a triplet $(M^n, g, V)$, corresponding to a spin asymptotically flat static vacuum with an outward-directed equipotential 
photon surface $\Sigma$ as boundary. Then, it follows from \cite[Formulas 4.16 and 4.17]{CederbaumGalloway2} that $\Sigma$ is a totally 
umbilical hypersurface with constant mean curvature, for which (\ref{ScalarMeanCurvatureConditionRig}) holds, and $V$ satisfies the overdetermined problem 
(\ref{OverBP}) with $\Gamma = c - 1$. Thus, Theorem \ref{RigidityGeneral} applies directly.

To place our result in context, let us briefly review the uniqueness results in this area. The first results in this direction focused on photon spheres and were obtained by Cederbaum \cite{Cederbaum1} and Cederbaum-Galloway \cite{CederbaumGalloway1} in the $(3+1)$-dimensional case. Subsequently, several generalizations have been obtained under various conditions. From the perspective of the inner boundary conditions, the most general version of this uniqueness result (to the author's knowledge) can be found in \cite{Cederbaum2, CederbaumGalloway2, Raulot12}, where the boundary may be disconnected and include black holes and outward-directed equipotential photon components. The proof relies on the positive mass theorem, which requires $3 \leq n \leq 7$ (though see \cite{Lohkamp2, SchoenYau3}) or that $M$ is spin, as is our case here. However, the asymptotic conditions on the metric and static potential are more restrictive, as asymptotic isotropy is assumed, while in
our result, only asymptotic flatness is required. More recently, Cederbaum, Cogo, Leandro, and Dos Santos \cite{CederbaumCogoLeandroDosSantos} generalized 
Robinson's method \cite{Robinson} (originally developed to prove uniqueness theorems for static black holes) to the case of photon surfaces. In their work, 
the assumptions on the asymptotic behavior of the metric and the static potential are much weaker than in all previous studies. This approach also allows 
one to dispense with the "outward-directed" assumption. However, the boundary is assumed to be connected, and a strong condition on the scalar curvature is 
required when $n \geq 4$. We should also mention the recent works of Harvie and Wang \cite{HarvieWang2,HarvieWang1}, where a new approach based on the static Minkowski inequality is introduced but where a stability condition and a condition on the Einstein-Hilbert of $\Sigma$ are required. The assumptions we impose in Theorem \ref{PhotonSurfacesUniqueness} are therefore standard in terms of asymptotic behavior and very general from the perspective of the boundary geometry (except for the connectedness).

\subsection{An uniqueness result for asymptotically flat static manifolds with boundary} In \cite{CruzVitorio}, Cruz and Vittorio introduced the notion of {\it static manifolds with boundary}. On a manifold with boundary, when studying the surjectivity of the map that associates to a Riemannian metric $g$ the couple $(R,H)$, where $R$ and $H$ are respectively the scalar and the mean curvatures of $g$, we are led to examine the kernel of the formal $L^2$-adjoint of its linearized operator. Then it can then be checked that a function $V$ is in this kernel if it satisfies the following boundary problem
$$.
\left\lbrace
\begin{array}{rll}
\nabla^2 V-(\Delta V)g-V{\rm Ric} & = & 0\text{ on }M \\
\frac{\partial V}{\partial \nu}\gamma-V A & = & 0 \text{ on }\Sigma.
\end{array}
\right.
$$ 
A manifold with boundary which possesses such a non-trivial function is called a {\it static manifold with boundary}. Recently, Medvedev \cite[Theorem 1.11]{Medvedev} proved an uniqueness result (among other things) for asymptotically isotropic static manifolds with equipotential boundary. Theorem \ref{RigidityGeneral} ensures that his result still holds if we only assume that the manifold is spin, asymptotically flat and the boundary is connected. 
\begin{theorem}
Let $(M^n,g,V)$ be a one-ended asymptotically flat spin static manifold with a compact and connected equipotential boundary and denote by $m$ its ADM mass. If $V$ is non-zero on $\partial M$ and $H>0$ then $(M^n,g)$ is isometric to the exterior of the (unique) photon sphere in the Schwarzschild manifold of mass $m>0$.
\end{theorem}
The exterior of the unique photon sphere is the part of the Schwarzschild manifold defined by $\mathbb{M}^n_m(r_S)$ with $r_S=\big(n-1+\sqrt{n(n-2)}\big)/2$. Once again, this result follows directly from Theorem \ref{RigidityGeneral} with $c=n/(n-2)$ (see \cite[Lemma 6.3]{Medvedev}). 


\section{Proof of Theorem \ref{SpinCase} and Theorem \ref{ThreeDimensionalCase}}\label{TheProofs}


The proofs of Theorems \ref{SpinCase} and Theorem \ref{ThreeDimensionalCase} relies on a mass-capacity inequality by Hirsch and Miao \cite{HirschMiao1} and on positive mass for manifolds with boundary (see Appendix \ref{AppendixA}). 

\vspace{0.2cm}

Here we set
\begin{eqnarray*}
0<\alpha:=\frac{\beta}{1+\beta}<1\quad\text{where}\quad\beta:=\frac{\Lambda}{c-1}>0,
\end{eqnarray*}
and so the condition (\ref{AdditionalBC}) now reads as 
\begin{equation}\label{AdditionalBC1}
-2\frac{1-\alpha}{\alpha}{\frac{\partial \Phi}{\partial\nu}}=(c-1)\frac{n-2}{n-1}H\quad\text{on }\Sigma.
\end{equation}

From \cite{HirschMiao1} we know that on an $n$-dimensional asymptotically flat manifold $(M^n,g)$ with nonnegative scalar curvature and with compact inner boundary, the mass-capacity inequality
$$
m\geq (1-\alpha)\,\mathcal{C}(\Sigma,M)
$$
holds for $0\leq\alpha\leq 1$ if
\begin{equation}\label{HM-Condition}
-\frac{2\alpha}{1+\alpha}\frac{\partial\Phi}{\partial\nu}\geq\frac{n-2}{n-1}H.
\end{equation}
Moreover, equality occurs if and only if $(M^n,g)$ is isometric to the exterior region outside a rotationally symmetric sphere in a Riemannian Schwarzschild manifold. 
Then, if we first assume that $\alpha^2c\geq 1$, it is straightforward to see that $c-1\geq (1-\alpha^2)\alpha^{-2}$ and so it follows from (\ref{AdditionalBC1}) that the condition (\ref{HM-Condition}) is fulfilled. This implies that our mass-capacity is true in this situation as well as the equality case. 

Let us now show that this also holds if $\alpha^2c\leq 1$. For this, consider the metric conformally related to $g$ defined by 
\begin{eqnarray*}
\overline{g}_\alpha=\Psi_\alpha^{4/(n-2)} g\quad\text{where}\quad\Psi_\alpha:=1-\frac{1-\alpha}{2}\Phi
\end{eqnarray*}
is a smooth positive function since $0<\Phi\leq 1$ and $\alpha<1$. From the asymptotic expansion (\ref{ExpansionBCP}) of $\Phi$, it is immediate to check that $(M^n,\overline{g}_\alpha)$ is an asymptotically flat manifold with ADM mass given by
\begin{eqnarray}\label{MassConformal}
\overline{m}_\alpha=m-(1-\alpha)\,\mathcal{C}(\Sigma,M).
\end{eqnarray}
Moreover, the scalar curvature of $(M^n,\overline{g})$ is easily computed to be
\begin{eqnarray*}
\overline{R}_\alpha=\Psi_\alpha^{-\frac{n+2}{n-2}}\Big(-4\frac{n-1}{n-2}\Delta\Psi_\alpha+R\Psi_\alpha\Big)\geq 0
\end{eqnarray*}
since $\Psi_\alpha$ is harmonic and $R$ is nonnegative. As well, the mean curvature of $\Sigma$ in $(M^n,\overline{g})$ is 
\begin{eqnarray*}
\overline{H}_\alpha=\Psi_\alpha^{-\frac{n}{n-2}}\Big(2\frac{n-1}{n-2}\frac{\partial \Psi_\alpha}{\partial\nu}+H\Psi_\alpha\Big)
\end{eqnarray*}
which can be rewritten as
\begin{eqnarray}\label{MeanCurvatureAlpha}
\overline{H}_\alpha=2^{\frac{2}{n-2}}(1+c\alpha)(1+\alpha)^{-\frac{n}{n-2}}H
\end{eqnarray}
since $\Psi_\alpha$ is constant equals to $(1+\alpha)/2$ on $\Sigma$ and since because of (\ref{AdditionalBC1}) it holds on $\Sigma$ that
\begin{eqnarray*}
\frac{\partial\Psi_\alpha}{\partial\nu}=\frac{\alpha(c-1)}{4}\frac{n-2}{n-1}H.
\end{eqnarray*}
On the other hand, observe that the metrics $\overline{g}_\alpha$ and $g$ being homothetic when restricted to $\Sigma$, it is easy to check that their corresponding scalar curvatures are related by the formula
\begin{eqnarray}\label{ScalarBoundaries}
\overline{\Ss}_\alpha=2^{\frac{2}{n-2}}(1+\alpha)^{-\frac{2}{n-2}}\Ss.
\end{eqnarray}
Notice that, until now, these arguments apply for the proof of both Theorems \ref{SpinCase} and \ref{ThreeDimensionalCase}. 

Let us now look at the former case. In fact, we assume that the assumption (\ref{ScalarMeanHypothese}) in Theorem \ref{SpinCase} holds,  and so we deduce that
\begin{eqnarray*}
\overline{\Ss}_{\alpha}\geq 2^{\frac{2}{n-2}}(1+\alpha)^{-\frac{2}{n-2}}\inf_\Sigma(\Ss)\geq 2^{\frac{2}{n-2}}(1+\alpha)^{-\frac{2}{n-2}}c\frac{n-2}{n-1}\max_\Sigma(H^2),
\end{eqnarray*} 
which from (\ref{MeanCurvatureAlpha}) rewrites as 
\begin{eqnarray*}
\overline{\Ss}_{\alpha}\geq\frac{c(1+\alpha)^2}{(1+c\alpha)^2}\left(\frac{n-2}{n-1}\right)\max_\Sigma(\overline{H}_\alpha^2).
\end{eqnarray*} 
However, since $\alpha^2c\leq 1$, it is immediate to see that $c(1+\alpha)^2\geq(1+c\alpha)^2$ so that 
\begin{eqnarray*}
\inf_\Sigma(\overline{\Ss}_\alpha)\geq\frac{n-2}{n-1}\max_\Sigma(\overline{H}^2_\alpha).
\end{eqnarray*} 
We thus have shown that the manifold $(M^n,\overline{g}_{\alpha})$ satisfies the assumptions of Theorem \ref{PMT-Boundary} in Appendix \ref{AppendixA} and then its ADM mass $\overline{m}_\alpha$ is nonnegative. The conclusion follows directly from (\ref{MassConformal}). Now if equality is achieved, the manifold $(M^n,\overline{g}_{\alpha})$ reaches the equality case in Theorem \ref{PMT-Boundary} and so it is isometric to the exterior of a round sphere with radius $r_\alpha>0$ in the Euclidean space. Then it follows that the smooth function $\Psi_{\alpha}^{-1}$ satisfies
$$
\left\lbrace
\begin{array}{ll}
\Delta_\delta\Psi^{-1}_\alpha= 0 & \text{ in } \mathbb{R}^n\setminus B(0,r_\alpha) \\
\Psi^{-1}_\alpha=2/(1+\alpha) & \text{ on } \mathbb{S}^{n-1}_{r_\alpha}\\ 
\Psi^{-1}_\alpha\rightarrow 1 & \text{ as } |x|\rightarrow\infty
\end{array}
\right.
$$
where $\Delta_\delta$ is the Laplace operator in the Euclidean space. It turns out that the unique solution of the aforementioned boundary problem is easily seen to be 
\begin{eqnarray*}
\Psi^{-1}_\alpha(x)=1+\frac{m}{2r^{n-2}}\quad\text{with}\quad m=2\,\frac{1-\alpha}{1+\alpha}r_\alpha^{n-2}
\end{eqnarray*}
and thus $(M^n,g)$ is isometric to $\big(\mathbb{M}^n_m(r_\alpha),g_m\big)$. Conversely, it follows from (\ref{MeanCurvatureSS}) and (\ref{SchwarzschildCapacitors}) that $\mathbb{M}^n_m(r_0)$, the region outside a rotationally symmetric sphere exterior in an $n$-dimensional Schwarzschild manifold of mass $m>0$, with $r_0>(m/2)^{1/(n-2)}$ is a critical area-normalized capacitor with 
\begin{eqnarray*}
 H_{m,r_0}>0\quad\text{and}\quad\Lambda_{m,r_0}=\frac{2}{1-\frac{m}{2r^{n-2}_0}}.
\end{eqnarray*}
Moreover, we compute from (\ref{CapacityPotentialSS}) that the capacity of $\Sigma_{r_0}$ in $\mathbb{M}^n_m(r_0)$ is 
\begin{eqnarray*}\label{SchwarzschildCapacity}
\mathcal{C}\big(\Sigma_{r_0},\mathbb{M}^n_m(r_0)\big)=\frac{m}{2}+r_0^{n-2}.
\end{eqnarray*} 
Finally from (\ref{RadiusSS}) and (\ref{MeanCurvatureSS}), we observe that if $\Ss_{m,r_0}$ denotes the scalar curvature of $\Sigma_{r_0}$ we have
\begin{eqnarray*}
\Ss_{m,r_0}=\frac{n-2}{n-1}c H_{m,r_0}^2\quad\text{with}\quad c:=\left(\frac{1+\frac{m}{2r_0^{n-2}}}{1-\frac{m}{2r_0^{n-2}}}\right)^2>1
\end{eqnarray*}  
and so equality occurs in (\ref{MassCapIneq}) for $\mathbb{M}^n_m(r_0)$. This conclude the proof of Theorem \ref{SpinCase}. 

\vspace{0.2cm}

We are now left to finish the proof of Theorem \ref{ThreeDimensionalCase} where we assumed (\ref{ScalarMeanHypothese1}) instead of (\ref{ScalarMeanHypothese}). In this case, the formulae (\ref{MeanCurvatureAlpha}) and (\ref{ScalarBoundaries}) ensure that $\overline{\Ss}_\alpha\geq\overline{H}^2_\alpha/2$. We thus have shown that the manifold $(M^n,\overline{g}_\alpha)$ satisfies the assumptions of Theorem \ref{PMT-Boundary2} in Appendix \ref{AppendixA} and then its ADM mass $\overline{m}_\alpha$ is nonnegative. The conclusion follows directly from (\ref{MassConformal}). The equality case follows exactly as above.


\appendix



\section{Two positive mass theorems on manifold with boundary}\label{AppendixA}


In this section, we give the proof of the two positive mass theorems for manifolds with boundary we used to demonstrate our main results. As we will see, they are direct consequences on previous works of Herzlich \cite{Herzlich1,Herzlich2} and Friedrich \cite{Friedrich3} for the spin case and of Miao \cite{Miao1} for the other cases.  

\vspace{0.2cm}

Let us first consider the spin case. In this situation, we have:
\begin{theorem}\label{PMT-Boundary}
Let $(M^n,g)$ be an asymptotically flat spin manifold with a compact and connected non empty inner boundary $\Sigma^{n-1}:=\partial M^n$. Assume that $R\geq 0$ and that $H>0$. If 
\begin{equation}\label{ScalarMeanHypothese_EuclideanSpace}
\inf_\Sigma(\Ss)\geq\frac{n-2}{n-1}\max_{\Sigma}(H^2)
\end{equation}
then $m\geq 0$ and equality occurs if and only if $(M^n,g)$ is isometric to the exterior of a round sphere in the Euclidean space. 
\end{theorem}

{\it Proof.} In \cite[Proposition 2.1]{Herzlich1} and \cite[Proposition 2.1]{Herzlich2}, Herzlich proved that the mass of an asymptotically flat manifold with compact boundary is nonnegative if the first eigenvalue of the Dirac operator of $\Sigma$, denoted by $\lambda_1$, satisfies the inequality $\lambda_1\geq\max_{\Sigma}(H)/2$. On the other hand, the Friedrich inequality \cite{Friedrich3} ensures that $\lambda_1$ satisfies the lower bound
\begin{eqnarray*}
\lambda_1^2\geq\frac{n-1}{4(n-2)}\,\inf_\Sigma(\Ss)
\end{eqnarray*}
and so the nonnegativity of the mass under the assumption (\ref{ScalarMeanHypothese_EuclideanSpace}) follows from the combination of these two results. Now if the mass is zero, we argue as in \cite[p.14]{Herzlich2}. First, from Herzlich's result, we deduce that $(M^n,g)$ is a flat manifold with constant mean curvature boundary. On the other hand, we also have equality in the Friedrich inequality and so $\Sigma$ carries a real Killing spinor. It is in particular an Einstein manifold with constant positive scalar curvature. Then if follows from the Gauss formula that $\Sigma$ has to be a totally umbilical hypersurface with constant mean curvature in $M$. Since $M$ is flat, the Gauss formula once again ensures that $\Sigma$ has constant sectional curvature and so it has to be a quotient of a round sphere. From the Heintze-Karcher inequality, we deduce that $\Sigma$ is in fact a round sphere. Then we can glue to $M$ along $\Sigma$ an Euclidean ball whose boundary is isometric to $\Sigma$ to obtain a smooth and complete asymptotically flat manifold with zero scalar curvature and zero mass for which the rigidity part of the classical positive mass theorem applies. This allows to conclude that $(M^n,g)$ is isometric to the exterior of a round ball in the Euclidean space.  
\hfill$\square$

\vspace{0.2cm}

One could get ride of the spin assumption by proving a positive mass theorem for manifolds with an isolated conical singularity and admitting corners along a hypersurface, that is a version {\it \`a la Miao} \cite[Theorem 1]{Miao1} of \cite[Theorem 1.1]{DaiSunWang}. Then one can proceed as follow: consider the product manifold $\Omega:=(0,r_0]\times\Sigma$ with $r_0=(n-1)(\max_\Sigma H)^{-1}$ endowed with the warped product metric $\widetilde{g}=dr^2+r^2 r_0^{-2}\gamma$ where $\gamma:=g_{|\Sigma}$. This Riemannian manifold $(\Omega,\widetilde{g})$ has a conical singularity in $\{0\}\times\Sigma$ and has the following properties:
\begin{enumerate}[label=(\roman*)]
\item the scalar curvature $\widetilde{R}$ of $\widetilde{g}$ is given by 
\begin{eqnarray*}
\widetilde{R}=r_0^2r^{-2}\left(\Ss-(n-1)(n-2)r_0^{-2}\right)
\end{eqnarray*} 
and satisfies $\widetilde{R}\geq 0$ because of (\ref{ScalarMeanHypothese_EuclideanSpace});
\item the metric induced by $\widetilde{g}$ on the slice $\{r_0\}\times\Sigma$ is $\gamma$;
\item the mean curvature of $\{r_0\}\times\Sigma$ in $\Omega$ is $\widetilde{H}=\max_\Sigma H$.
\end{enumerate}
In this way, we obtain a fill-in with an isolated conical singularity which can be glued to $M$ along $\Sigma$ to get an asymptotically flat manifold with a corner along $\Sigma$ satisfying $\widetilde{H}\geq H$ and with an isolated conical singularity. Moreover, the scalar curvature of this new manifold is nonnegative on the smooth part so that this positive mass theorem would apply.  

\vspace{0.2cm}

An alternative boundary condition is provided below to ensure nonnegativity of the mass.
\begin{theorem}\label{PMT-Boundary2}
Let $(M^n,g)$ be an asymptotically flat manifold with a compact and connected non empty inner boundary $\Sigma^{n-1}:=\partial M^n$. Assume also that $R\geq 0$, that $H>0$ and that 
\begin{equation}\label{ScalarMeanHypothese_EuclideanSpace1}
\Ss\geq\frac{n-2}{n-1} H^2
\end{equation}
where $\Ss$ is the scalar curvature of $\Sigma$. If one of the following conditions is fulfilled:
\begin{enumerate}
\item $n=3$;
\item $(\Sigma^{n-1},\gamma)$ admits an isometric embedding in the flat $n$-dimensional Euclidean space;
\end{enumerate}
then $m\geq 0$ and equality occurs if and only if $(M^n,g)$ is isometric to the exterior of a round sphere in the Euclidean space. 
\end{theorem}

{\it Proof.} First note that when $n=3$, since $M$ is assumed to be orientable so is $\Sigma$ and from the assumption (\ref{ScalarMeanHypothese_EuclideanSpace1}), we deduce that $\Sigma$ is a topological $2$-sphere with positive Gauss curvature. It follows from the Weyl's embedding theorem \cite{Nirenberg,Pogorelov} that $(\Sigma,\gamma)$ admits an unique strictly convex isometric embedding in the Euclidean space $(\mathbb{R}^{3},\delta)$. This means that the assumption {\it (2)} is automatically satisfied for $n=3$. So for $n\geq 3$, the Gauss formula for $\Sigma$ with respect to its Euclidean embedding ensures that 
\begin{eqnarray}\label{GaussEuclidean}
\Ss=\frac{n-2}{n-1}H_0^2-|\mathcal{O}_0|^2\leq \frac{n-2}{n-1}H_0^2
\end{eqnarray}
where $\mathcal{O}_0:=A_0-\frac{H_0}{n-1}\gamma$ is the traceless part of the second fundamental form of $\Sigma$ in $\mathbb{R}^n$. Here $A_0$ and $H_0$ are respectively the associated second fundamental form and mean curvature. In particular, from (\ref{ScalarMeanHypothese_EuclideanSpace1}), we get that $|H_0|\geq H$ since $H>0$. However, since $\Sigma$ is compact, there is point $p\in\Sigma$ such that $H_0(p)>0$ so that $H_0\geq H$ holds everywhere on $\Sigma$. On the other hand, since $\Sigma$ is compact and embeds in $\mathbb{R}^n$, it bounds a compact domain say $\Omega$. Then we glue $\Omega$ to $M$ along $\Sigma$ to get an asymptotically flat manifold with nonnegative scalar curvature on its smooth part and a corner along $\Sigma$ with $H_0\geq H$ so that Miao's positive mass theorem \cite[Theorem 1]{Miao1} allows to conclude that the mass of $(M^n,g)$ is nonnegative. Finally, if the mass is zero, it follows from \cite{McFeronSzekelyhidi} that $(M^n,g)$ is isometric to the unbounded domain delimited by $\Omega$ in $\mathbb{R}^n$. Moreover, it also holds that $H_0=H$ and so equality occurs in (\ref{GaussEuclidean}). Then $\Sigma$ has to be totally umbilical in the Euclidean space and so isometric to a round sphere. This concludes the proof. 
\hfill$\square$


\bibliographystyle{alpha}     
\bibliography{BiblioHabilitation}


\end{document}